\documentclass[10pt,twocolumn,twoside]{IEEEtran}
\usepackage[cmex10]{amsmath}

\usepackage{multirow}
\usepackage{romannum}
\usepackage{nicematrix}
\usepackage{cite}
\usepackage{amssymb,amsfonts}
\usepackage{amsthm}
\usepackage{bm}
\usepackage{algorithmic}
\usepackage{graphicx}
\usepackage{textcomp}
\usepackage{xcolor}
\usepackage{mathabx}
\usepackage{subfiles} 
\usepackage{caption} 
\usepackage{enumitem}
\usepackage{subcaption}
\usepackage{nccmath, mathtools}
\usepackage{url}

\newcommand{\st}{\text{s.t.}}

\begin{document}

\newtheorem{definition}{Definition}

\newtheorem{lemma}{Lemma}

\newtheorem{theorem}{Theorem}

%
\title{Convex Restriction of Feasible Sets for AC Radial Networks}

\author{
\IEEEauthorblockN{Ling Zhang, Daniel Tabas, and Baosen Zhang\\ }
\IEEEauthorblockA{Electrical and Computer Engineering, University of Washington\\
\{lzhang18, dtabas, zhangbao\}@uw.edu}
}


\maketitle

Many problems in power systems involve optimizing a certain objective function subject to power flow equations and engineering constraints.
A long-standing challenge in solving them is the nonconvexity of their feasible sets.
In this paper, we propose an analytical method to construct the convex restriction 
of the feasible set for AC power flows in radial networks. The construction relies on simple geometrical ideas and is explicit, in the sense that it does not involve solving other complicated optimization problems. We also show that the construct restrictions are in some sense maximal, that is, the best possible ones. 
Optimization problems constrained to these sets are not only simpler to solve but also offer feasibility guarantee for the solutions to the original OPF problem.
Furthermore, we present an iterative algorithm to improve on the solution quality by successively constructing a sequence of convex restricted sets and solving the optimization on them.
The numerical experiments on the IEEE 123-bus distribution network
show that our method finds good feasible solutions within just a few iterations and works well with various objective functions, even in situations where traditional methods fail to return a solution.



\section{Introduction}\label{sec:intro}


Many problems in power systems determine an optimal
network operating point that seeks to minimize a certain objective, all while satisfying a set of power flow equations and engineering constraints, 
such as optimal power flow (OPF), state estimation, voltage regularization, etc. 
However, the inherent nonlinearity of power flow equations gives rise to nonconvex feasible sets for these problems, which makes even finding feasible solutions nontrivial~\cite{lehmann2016NPhard}.

One strategy to address the nonconvexity of OPF problems is to convexify the constraints. This results in convex problems that are simpler to solve. There are mainly two convexification approaches.
First is convex relaxation,  which finds an outer approximation of the feasible set by modeling the original problem as a semidefinite or a conic program~\cite{jabr2006scop,bai2008sdp,low2014, low2014part2}. If the solution obtained from relaxation is also a feasible solution for the original problem, then it is a globally optimal solution~\cite{lavaei2011zero,lavaei2014geometry}. This implies that the relaxation is tight.
However, in practice, the convex relaxation can lead to non-physical solutions when it is not tight (this is common when bus power lower bounds are binding)~\cite{zhang2013geometry}. In such situations, distinguishing between whether the original problem is genuinely infeasible or if the relaxation method has failed becomes challenging.

Unlike convex relaxation, the second approach, convex restriction, provides an inner approximation to the feasible region. 
Optimization within these convex subsets guarantees
the feasibility of solutions for the original OPF problem. In essence, if the convex restriction algorithm produces a solution, that solution is guaranteed to be physically attainable~\cite{wu1982fp,simpson-porco2017fppf}.
In this paper, we focus on convex restrictions of OPF feasible sets in radial networks.
It's noteworthy that most distribution networks are operated radially.
Due to the growing need for integrating distributed generation and facilitating demand response~\cite{ipakchi2009futuregrid}, solving OPF in distribution networks has become increasingly important.

Most existing convex restriction approaches work in the power injection space~\cite{nguyen2019inner,lee2019restriction}. 
It turns out working in this space often requires an assumption that the admissible power injection space or voltage space (or both) has a polytopic shape for analytical convenience.
Determining a non-conservative polytopic inner approximation of the original feasible set can be nontrivial. It may involve solving nonconvex optimization problems~\cite{nguyen2019inner}. Attempts to simplify or bypass this computationally intensive step often result in overly conservative results (see e.g., Figure~\ref{fig:3-bus} as an example). Moreover, this assumption itself may be overly restrictive as it approximates every nonlinear constraint using a linear one. 

To address these challenges, we propose to construct the convex restriction in a transformed coordinate space of voltage phase angles. Specifically, we apply a change of variables such that the active and reactive power equations become naturally convex after variable change, hence eliminating the need to approximate the bus power upper bound constraints. The lower bound constraints can be approximated using the first-order Taylor approximation, which is the best (the least conservative) upper bound one could have for concave functions. We use a 3-bus line network as an example (Figure~\ref{fig:3-bus}) to show that the convex restriction constructed in this way can be a maximal convex subset, meaning it cannot be contained within any other convex subset.

One may note that different tangent points can produce different first-order Taylor approximations, thus leading to different convex restricted sets. Therefore, we introduce an iterative algorithm that progressively refines the locations of tangent points to find the optimal convex restriction, one that contains the optimal solution of the original problem.
We test our method on the IEEE 123-bus distribution network with three types of objective functions applied: power loss minimization, generation cost minimization, and state estimation. The simulation results show that the iterative algorithm that builds on the proposed convex restriction method 
always finds a good feasible solution within at most $10$ iterations, even when the traditional methods failed.

The paper is organized as follows. In Section II, we present the original formulation of the AC optimal power flow problem.
Section III introduces our proposed convex restriction method and the resulting convex restricted OPF. We examine the geometry of our constructed convex restricted set in Section IV. Based on the insights from Section IV, we present the iterative algorithm in Section V to obtain an optimal solution for the original OPF problem.
Additionally, in Section V, we also compare the proposed method against several baseline methods using the IEEE 123-bus distribution network. Finally, Section VI concludes the paper.

\section{Optimal Power Flow}\label{sec:settings}
Consider a radial network, where
$\mathcal{N}$ is the set of buses and $\mathcal{E}$ the set of lines. For simplicity, we assume that the voltage magnitudes are at 1 p.u. and consider the following problem:
\begin{subequations}\label{eqs:opf}
    \begin{align}
     \min_{\bm{\theta}} \, & c(\mathbf{p},\mathbf{q})\label{opf:objective}\\
     \st \,
     & p_i = \textstyle \sum_{j:(i,j)\in\mathcal{E}} g_{ij} - g_{ij}\cos{(\theta_{ij})} 
    +  b_{ij}\sin{(\theta_{ij})}\label{opf:p}\\
    & q_i = \textstyle \sum_{j:(i,j)\in\mathcal{E}} b_{ij} - b_{ij}\cos{(\theta_{ij})} 
    -  g_{ij}\sin{(\theta_{ij})}\label{opf:q}\\
     & \underline{P}_i \leq p_i\leq \overline{P}_i, \; \underline{Q}_i \leq q_i\leq \overline{Q}_i, \; \underline{\theta}_{ij}\leq \theta_{ij}\leq \overline{\theta}_{ij}\label{opf:cons}    \end{align}   
\end{subequations}
where $\bm{\theta}$ is the voltage angle vector, $\theta_{ij}=\theta_i-\theta_j$ is the angle difference between bus $i$ and $j$, $g_{ij} -j b_{ij}$  is the admittance of the line $(i,j)$,   and $p_i$ and $q_i$ are active and reactive power at each bus $i$, respectively. We do not explicitly specify whether a bus is a generator or a load. This information can be inferred from the upper and lower bounds on power, for example, if $\overline{P}_i$ is negative, then bus $i$ is a load bus. We will come back to the objective in later sections, but it suffices to think of $c$ as some cost function in active and reactive power. 

The constraint $\underline{\theta}_{ij}\leq \theta_{ij}\leq \overline{\theta}_{ij}$ is important to our analysis. Particularly, we assume that the limits satisfy $[\underline{\theta}_{ij},\overline{\theta}_{ij}] \subset (-\frac{\pi}{2},\frac{\pi}{2})$. Under this assumption, $\sin(\theta_{ij})$ is a monotonic function of $\theta_{ij}$, which forms the basis of our approach. We believe this assumption is likely to be true in practice, since it is difficult to think of a distribution system where the angle differences would be larger than 90 degrees. 

We denote the feasible set of \eqref{eqs:opf} as $\Theta$. Since $\Theta$ is not convex in general, solving \eqref{eqs:opf} is nontrivial and a number of numerical methods have been developed~\cite{sun1984optimal,lavaei2011zero,zhang2012geometry,low2014convex,molzahn2019survey}. A drawback of all these methods is that if they do not return a solution--for example, when a convex relaxation is not tight or when a Newton-type algorithm fails to converge--it's not easy to tell whether the problem itself is infeasible or it's the algorithm that has failed.
Unlike these methods, convex restriction finds an inner subset of $\Theta$. If this convex subset is not empty, it serves as evidence that the original problem is indeed feasible, and optimization over this subset is guaranteed to produce a feasible solution.
In the next section, we introduce a simple change of variables technique that helps to construct a convex restriction of $\Theta$.

\section{Convex Restriction} \label{sec:restriction}
We introduce the following change of variables:
Let $z_{ij}=\sin(\theta_{ij})$, for all $(i,j)\in\mathcal{E}$. 
Note that this transformation preserves the feasibility and optimality of \eqref{eqs:opf}, namely, the problem in the transformed coordinates is  equivalent to \eqref{eqs:opf}.
Since the sine function is monotonically increasing in $(-\frac{\pi}{2}, \frac{\pi}{2})$ and is invertible, the angle difference $\theta_{ij}$ and hence the angles themselves $\bm{\theta}$ can be readily obtained from $z_{ij}$. Therefore, in the rest of the paper, we focus on solving the problem in the transformed coordinates.
Next, we look at how each of the constraints in \eqref{eqs:opf} are represented after this change of variables. 

\subsection{Angle Difference Constraints}
The angle difference constraints $\underline{\theta}_{ij}\leq \theta_{ij}\leq \overline{\theta}_{ij}$ become $\underline{\theta}_{ij}\leq \sin^{-1}(z_{ij})\leq \overline{\theta}_{ij}$. This appears to be nonconvex in $z_{ij}$, but a simple observation is that because $\sin$ is monotonically increasing,  the constraint is equivalent to $\underline{z}_{ij} \leq z_{ij} \leq \overline{z}_{j}$, with $\underline{z}_{ij} = \sin^{-1} \underline{\theta}_{ij}$ and $\overline{z}_{ij} = \sin^{-1} \overline{\theta}_{ij}$, which are linear inequalities (and hence convex). 

\subsection{Active and Reactive Power Constraints}
Using the simple fact that for $\theta \in (-\frac{\pi}{2},\frac{\pi}{2})$, $\cos(\theta)=\sqrt{1-\sin(\theta)^2}$, the power flow from bus $i$ to bus $j$ is
$$ p_{ij}(z_{ij})= g_{ij} - g_{ij}\sqrt{1-z_{ij}^2} +  b_{ij}z_{ij} $$
in the $z_{ij}$ variables. Then letting $\mathbf{z}=\{z_{ij},(i,j) \in \mathcal{E}\}$, the nodal active power injection at bus $i$ is
\begin{equation} \label{eqn:active_z}
p_i(\mathbf{z}) = \sum_{j:(i,j)\in\mathcal{E}} g_{ij} - g_{ij}\sqrt{1-z_{ij}^2}  +  b_{ij}z_{ij}, 
\end{equation}
with the relevant constraint being $\underline{P}_i \leq p_i (\mathbf{z}) \leq \overline{P}_i$. 

The affine terms in \eqref{eqn:active_z} cause no difficulty. The nonlinear term, $-g_{ij}\sqrt{1-z_{ij}^2}$, is more interesting. By elementary calculations, $\sqrt{1-z_{ij}^2}$ is a concave function of $z_{ij}$.  Therefore, $- g_{ij}\sqrt{1-z_{ij}^2}$ is convex in $z_{i,j}$. Consequently, the upper bound on the active power, $p_i(\mathbf{z}) \leq \overline{P}_i$, is a convex constraint. In contrast to other  methods where every nonlinear constraint need to be approximated~\cite{nguyen2019inner,lee2019restriction} and causes the convex restriction to shrink, the upper bounds on active power are naturally convex in the $\mathbf{z}$ variables. This actually recovers a known result in OPF, where the problem tends to be convex under a condition called load over-satisfaction, meaning that the lower bounds are removed~\cite{baldick2006applied}.

The lower bound on active power has the form of a convex function greater than a constant, and is nonconvex. We replace it by a supporting hyperplane of the convex function, to create an upper estimate of the lower bound. Specifically, given a differentiable function $f$, it is convex if and only if the following first order condition is satisfied: 
\begin{equation} \label{eqn:convexity}
f(\mathbf{x}) \geq f(\mathbf{y})+\nabla f(\mathbf{y})^T (\mathbf{x}-\mathbf{y})
\end{equation}
for all $\mathbf{x}$ and $\mathbf{y}$ in its domain. Applying \eqref{eqn:convexity} to $p_i(\mathbf{z})$, we define:
\begin{equation} \label{eqn:p_i_lower}
 \underline{p}_i(\mathbf{z}) := p_i(\tilde{\mathbf{z}}_i^p)+\nabla p_i(\tilde{\mathbf{z}}_i^p)^T(\mathbf{z}-\tilde{\mathbf{z}}_i^p), 
\end{equation}
for some base point  $\tilde{\mathbf{z}}_i^p$. 
With \eqref{eqn:p_i_lower}, the two convex inequalities, $p_i(\mathbf{z}) \leq \overline{P}_i$ and $\underline{p}_i(\mathbf{z}) \geq \underline{P}_i$, together imply the original nonconvex inequality $\underline{P}_i \leq p_i(\mathbf{z}) \leq \overline{P}_i$. Particularly, in order to find a supporting hyperplane for each lower bound constraint, one must choose base points where these constraints are active. 
A valuable insight is that the constraint $p_i(\mathbf{z}) \leq \underline{P}_i$ is indeed convex. Therefore, 
a base point for each lower bound constraint can be found by starting with a strictly feasible point $\mathbf{z}^{\mathrm{o}}$ and projecting $\mathbf{z}^{\mathrm{o}}$ onto the lower bound constraints. Details of this procedure are described in the next section. 

The reactive power constrains can be treated in exactly the same way, by noticing that 
\begin{equation} \label{eqn:reactive_z}
q_i(\mathbf{z}) = \sum_{j:(i,j)\in\mathcal{E}} b_{ij} - b_{ij}\sqrt{1-z_{ij}^2}  -  g_{ij}z_{ij}, 
\end{equation}
is convex. The lower bounds can be handled by defining 
\begin{equation} \label{eqn:q_i_lower}
 \underline{q}_i(\mathbf{z}) := q_i(\tilde{\mathbf{z}}_i^q)+\nabla q_i(\tilde{\mathbf{z}}_i^q)^T(\mathbf{z}-\tilde{\mathbf{z}}_i^q),
\end{equation}    
at some base point $\tilde{\mathbf{z}}_i^q$. Then the convex restriction of the nonconvex inequality $\underline{Q}_i \leq q_i(\mathbf{z}) \leq \overline{Q}_i$ can be represented as two convex inequalities, $q_i(\mathbf{z}) \leq \overline{Q}_i$ and $\underline{q}_i(\mathbf{z}) \geq \underline{Q}_i$.

\subsection{OPF with a Convex Feasible Set}
All together, the convex restricted version of the problem in \eqref{eqs:opf} is
\begin{subequations} \label{eqs:convex_opf}
    \begin{align}
     \min_{\mathbf{z}} \, & c\big(\mathbf{p}(\mathbf{z}),\mathbf{q}(\mathbf{z})\big)\\
     \st \,
     & \eqref{eqn:active_z}, \ \eqref{eqn:reactive_z}, \label{constr:pz&qz}\\ 
     & \underline{\mathbf{z}} \leq \mathbf{z} \leq \overline{\mathbf{z}}, \label{constr:z_bounds}\\
     & p_i(\mathbf{z}) \leq \overline{P}_i, \ q_i(\mathbf{z}) \leq \overline{Q}_i,\label{constr:pz_bounds}\\
     & \underline{p}_i(\mathbf{z}) \geq \underline{P}_i, \ \underline{q}_i(\mathbf{z}) \geq \underline{Q}_i.\label{constr:qz_bounds}
     \end{align}   
\end{subequations}
The following theorem summarizes the main result of this section:
\begin{theorem} \label{thm:restriction}
   The feasible set of \eqref{eqs:convex_opf} is convex, and a feasible solution of \eqref{eqs:convex_opf} is a feasible solution of \eqref{eqs:opf}.
\end{theorem}

\captionsetup[figure]{font=small,skip=2pt}
\begin{figure}[ht!]
\centering
    \begin{subfigure}{0.3\textwidth}
    \includegraphics[width=6cm,height=6cm]{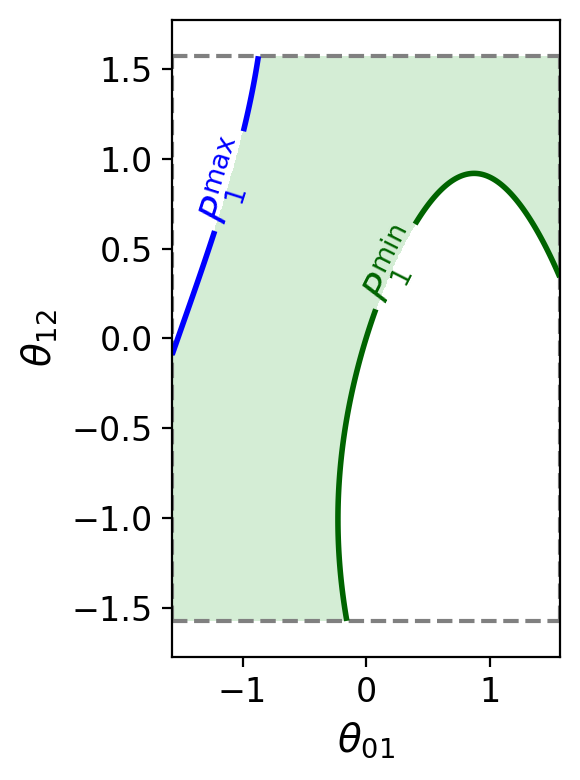}
    \label{fig:voltage_space}
    \end{subfigure}
    \begin{subfigure}{0.3\textwidth}
    \includegraphics[width=6cm,height=6cm]{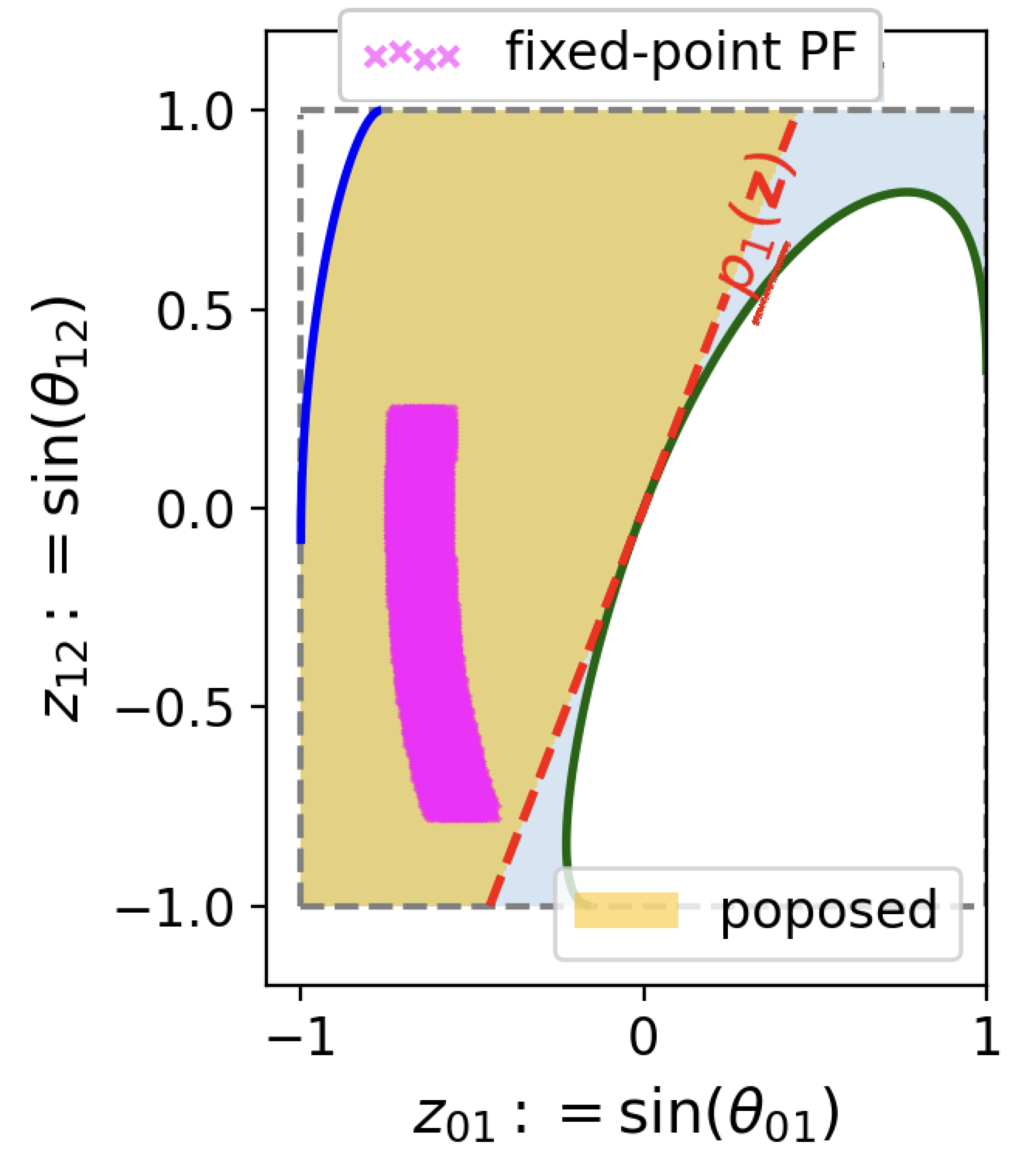}
    \label{fig:trans_space}
    \end{subfigure}
\caption{
A 3-bus line network example. Top is the feasible set in the $\theta$ space, and bottom is the feasible set in the $z$ space. Yellow is the convex restriction found by our method and magenta is the set found by the method in~\cite{lee2019restriction}. }
\label{fig:3-bus}
\end{figure}
In Figure \ref{fig:3-bus}, we use a three bus line network as an example to illustrate the discussions in this section. The top figure shows its 
feasible set in the $\theta$ space (green-colored area) and the bottom figure shows the transformed feasible set in the $z$ space, which is the cyan-colored area overlaid by the yellow-colored area. Particularly,
the yellow region in the bottom figure is the convex restriction constructed using the techniques introduced in this section. For comparison, the magenta-colored region, obtained from the method in~\cite{lee2019restriction} (it appears nonconvex because of the change in coordinates), is much more restricted than ours.

In fact, in this example, it is not possible to find another convex subset of the feasible region that contains the yellow region. Therefore, the yellow region can be considered as the maximal convex subset.
In the next section, we will take a closer look at the geometry of the feasible set for \eqref{eqs:convex_opf} and show why it is in some sense an optimal convex restriction.

\section{Geometry of Convex Restrictions}\label{sec:geometry}
We start with the following definition:
\begin{definition}[Maximal convex restriction] \label{def:maximal}
Let $f$ be a function from $\mathbb{R}^n$ to $\mathbb{R}$. Let $\mathcal{A}=\{\mathbf{x}: f(\mathbf{x}) \leq 0\}$ be a set in $\mathbb{R}^n$. We say a convex set $\mathcal{B}$ is a maximal convex restriction of $\mathcal{A}$ if 1) $\mathcal{B} \subseteq \mathcal{A}$,  and 2) there does not exist another convex set $\mathcal{C}$ such that $\mathcal{B} \subset \mathcal{C} \subset \mathcal{A}$.  
\end{definition}
Note that the bus power lower bound constraints all take the form of $\mathcal{A}=\{\mathbf{x}: f(\mathbf{x}) \leq 0\}$.
Definition \ref{def:maximal} says that given a set defined in this form, a convex restriction is maximal if it is not included in a larger convex subset of the set $\mathcal{A}$. The convex restriction shown in Fig.~\ref{fig:3-bus} illustrates this point. The yellow region at the bottom is maximal, since there does not exist another convex subset of the feasible region that contains it. But it is not the only maximal set, since the line that forms $\underline{p}_i$ could be tangent to the curve $p_i(\mathbf{z})=\underline{p}_i$ at other points. In fact, by changing the tangent lines, we obtain the family of maximal convex restrictions to the feasible set. We formalize this observation in the next theorem. 

\begin{theorem}\label{thm:maximal}
    If \eqref{eqs:opf} is strictly feasible, then each of the constraints defining the feasible set of \eqref{eqs:convex_opf} 
    forms a convex restriction 
    that is maximal in the sense of Definition \ref{def:maximal}. 
\end{theorem}
We want to emphasize that the maximal property is defined at a per-constraint level. Because multiple restricted sets can be maximal, the size of the feasible set in \eqref{eqs:convex_opf} depends on the choice of the base points. After stating the proof of Theorem~\ref{thm:maximal}, we outline an iterative procedure to select the base point that jointly optimizes the convex restrictions of all constraints. 

The proof of the theorem follows from the fact that the convex restriction of the bus power lower bound constraints is to find the best concave lower bound of a convex function,
which
is one of its supporting hyperplanes. More precisely, we have the following lemma:
\begin{lemma}~\label{lem:support}
    Let $f(\mathbf{x})$ be a convex function. Let $g(\mathbf{x})$ be a concave function. If $g(\mathbf{x}) \leq f(\mathbf{x})$ for all $\mathbf{x}$, then there exist an affine function $h(\mathbf{x})$ such that $g(\mathbf{x}) \leq h(\mathbf{x}) \leq f(\mathbf{x})$. 
\end{lemma}
\begin{proof}
First consider where the inequality is strict, that is, $g(\mathbf{x}) < f(\mathbf{x})$. A fundamental result in convex geometry is that there exist a separating hyperplane between $f$ and $g$~\cite{boyd2004convex}, and it serves as the function $h$. If the inequality is not strict, then pick a $\mathbf{y}$ such that $g(\mathbf{y})=f(\mathbf{y})$. Because $f$ is convex and $g$ is concave, there is a supporting hyperplane through the point $\mathbf{y}$ in the form of an affine function $h(\mathbf{x})$ such that $g(\mathbf{x}) \leq h(\mathbf{x}) \leq f(\mathbf{x})$. 
\end{proof}

Lemma \ref{lem:support} states that using affine functions to replace the active and reactive power lower bounds is the best one could do, and there does not exist another type of functions that will convexify the constraints while enlarging the feasibility region. Because the affine functions $\underline{p}_i(\mathbf{z})$ and $\underline{q}_i(\mathbf{z})$ are tangent to a point on the original lower bound, they are then maximal in the sense of Definition~\ref{def:maximal}. 

Next, we show how these tangent points can be found due to the following lemma, which requires the strict feasibility condition in the statement of Theorem~\ref{thm:maximal}. 
\begin{lemma}\label{lemma:projection}
    Let $\mathbf{z}^{\mathrm{o}}$ be a strictly feasible solution for \eqref{eqs:opf} and define
    $\underline{\mathcal{P}}_i = \{\mathbf{z}: p_i(\mathbf{z}) \leq \underline{P}_i\}$. 
    Let $\tilde{\mathbf{z}}_i^p$ be the Euclidean (2-norm) projection of $\mathbf{z}^{\mathrm{o}}$ onto $\underline{\mathcal{P}}_i$, i.e., $\tilde{\mathbf{z}}_i^p=\arg\min_{\mathbf{z}\in\underline{\mathcal{P}}_i} \|\mathbf{z}-\mathbf{z}^{\mathrm{o}}\|_2^2$. Then the projection $\tilde{\mathbf{z}}_i^p$
    lies on the curve $p_i(\mathbf{z})=\underline{P}_i$ whenever $\underline{\mathcal{P}}_i$ is nonempty. 
    
    By defining $\underline{\mathcal{Q}}_i = \{\mathbf{z}: q_i(\mathbf{z}) \leq \underline{Q}_i\}$, we have the same argument for reactive power equations that the projection onto $\mathcal{Q}_i$ lies on the curve $q_i(\mathbf{z})=\underline{Q}_i$.
\end{lemma}

\begin{proof}
    Since $\mathbf{z}^{\mathrm{o}}$ is a strictly feasible point, $\mathbf{z}^{\mathrm{o}} \not \in \underline{\mathcal{P}}_i$. Further, the set $\underline{\mathcal{P}}_i$ is convex. Due to these facts, the Euclidean projection of $\mathbf{z}^{\mathrm{o}}$ onto $\underline{\mathcal{P}}_i$ lies on the boundary of $\underline{\mathcal{P}}_i$. 
    {Now suppose for the sake of contradiction that the projection is
    an interior point of $\underline{\mathcal{P}}_i$, i.e., $p_i(\tilde{\mathbf{z}}_i^p) < \underline{P}_i$. Since the function $\hat{p}_i(\alpha) := p_i(\alpha \tilde{\mathbf{z}}_i^p + (1-\alpha)\mathbf{z}^{\mathrm{o}})$ is continuous for $\alpha \in [0,1]$, and  $\hat{p}_i(0) > \underline{P}_i$ and $\hat{p}_i(1) < \underline{P}_i$, there exists $\alpha^{\star} \in (0,1)$ such that $\hat{p}_i(\alpha^{\star}) = \underline{P}_i$. Let $\alpha^{\star} \tilde{\mathbf{z}}_i^p + (1-\alpha^{\star})\mathbf{z}^{\mathrm{o}}=\mathbf{z}^{\text{b}}$, and we have $\|\mathbf{z}^{\text{b}}-\mathbf{z}^{\mathrm{o}}\|_2^2=\|\alpha^{\star} \tilde{\mathbf{z}}_i^p + (1-\alpha^{\star})\mathbf{z}^{\mathrm{o}}-\mathbf{z}^{\mathrm{o}}\|_2^2 ={\alpha^{\star}}^2\|\tilde{\mathbf{z}}_i^p-\mathbf{z}^{\mathrm{o}}\|_2^2<\|\tilde{\mathbf{z}}_i^p-\mathbf{z}^{\mathrm{o}}\|_2^2$. This contradicts with $\tilde{\mathbf{z}}_i^p$ being the projection.} 
    The same proof logic applies to the reactive power equations.
\end{proof}

\begin{figure}[t]
\centering
 \begin{subfigure}{0.3\textwidth}
    \includegraphics[width=6.cm,height=6cm]{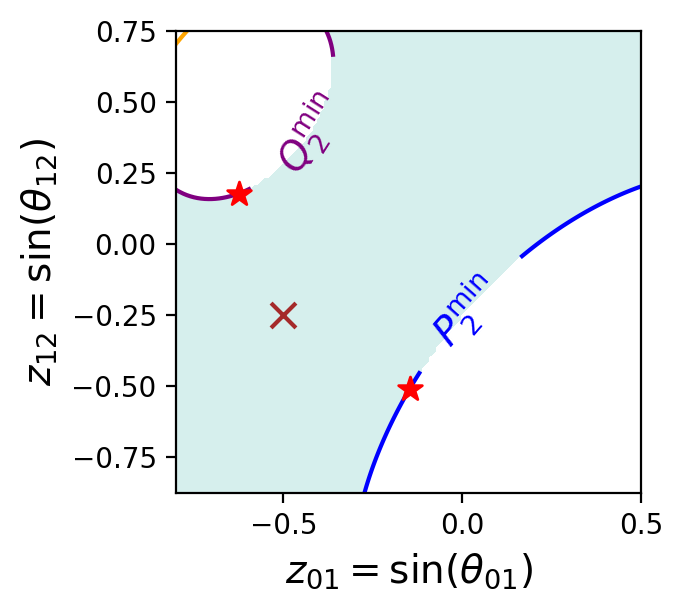}
    \label{fig:projection}
    \end{subfigure}
    \begin{subfigure}{0.3\textwidth}
    \includegraphics[width=6.cm,height=6cm]{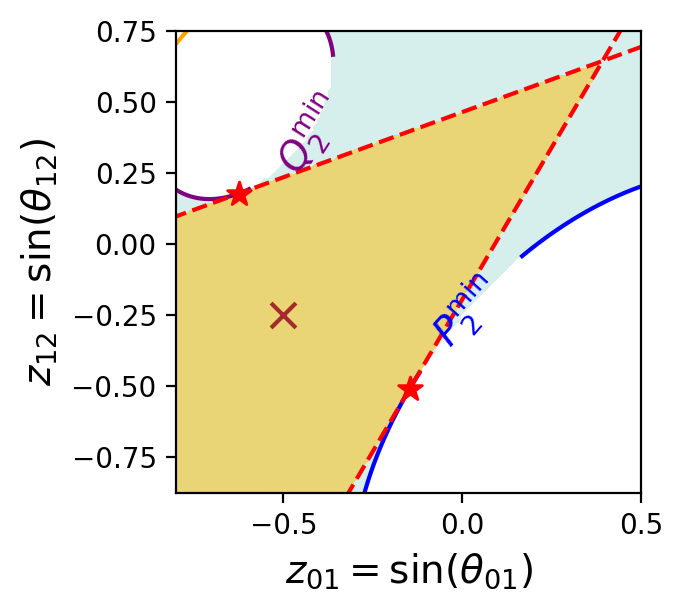}
    \label{fig:linearization}
    \end{subfigure}
\caption{
The cyan-colored region is the feasible set for a 3-bus line network and
the brown cross in it is a strictly feasible point. The red stars are its projections onto the convex sets $\underline{\mathcal{P}}_2 = \{\mathbf{z}: p_2(\mathbf{z}) \leq{P}_2^{\min}\}$ and $\underline{\mathcal{Q}}_2 = \{\mathbf{z}: q_2(\mathbf{z}) \leq{Q}_2^{\min}\}$, respectively.
In the bottom figure, the red dashed lines illustrate the linearization around these projected points and the resulting convex restriction is represented by the yellow region.
\vspace{-0.5cm}}
\label{fig:project-linearize}
\end{figure}

Lemma \ref{lemma:projection} suggests that we can obtain a series of tangent points by projecting a strictly feasible point onto each of the lower bound constraints in \eqref{constr:pz_bounds} and \eqref{constr:qz_bounds}.
Figure \ref{fig:project-linearize} shows the linearization of the lower bound constraints on active and reactive power for a 3-bus line network using the tangent points found by projection.
Note that if $\underline{\mathcal{P}}_i$ is empty, then it does not need to be linearized since the corresponding constraint $p_i(\mathbf{z}) \geq \underline{P}_i$ is vacuous and never active. This is easily checked when solving the projection. If the projection problem is infeasible, then a base point is not needed and the constraint can be removed. 

\captionsetup[figure]{font=small,skip=2pt}

It's important to note that by choosing a different initial feasible point $\mathbf{z}^{\mathrm{o}}$, a different set of tangent points can be obtained, hence resulting in a different convex restricted set. Building on this insight, in
the next section, we will introduce an iterative algorithm that progressively refines the choice of the feasible point and hence locations of the tangent points. As a result, a sequence of convex restricted sets are generated and each contains solutions progressively approach the optimal ones. We will evaluate the algorithm's performance using the IEEE 123-bus distribution network, considering various types of objective functions.

\section{Algorithm and Numerical Experiments}\label{sec:experiments}

In this section, we present the algorithm for solving 
\eqref{eqs:opf} by iteratively solving a sequence of \eqref{eqs:convex_opf}. In each iteration, we use the solution from the previous step to update the convex restricted feasible set for the next problem in the sequence.
Since these problems have convex feasible sets, finding a solution is always possible, which ensures the viability of this algorithm. The iterative procedure can be terminated once the obtained solutions stop changing or after a predetermined number of steps.

Note that the convexity of the overall problem in~\eqref{eqs:convex_opf} (both constraints and objective function) depends on the specific form of the objective function. We will show later that for a commonly used class of objective functions--positive linear combination of the generation costs--\eqref{eqs:convex_opf} is convex. For other cost functions, we will show that the iterative algorithm still performs very well. 

\subsection{Iterative Algorithm for Solving \eqref{eqs:opf}}
Starting from an initial feasible point $\mathbf{z}^{\mathrm{o},0}$,
denote the feasible point at the $k$-th iteration as $\mathbf{z}^{\mathrm{o},k}$.
Let $\{\tilde{\mathbf{z}}_i^{p,k}, \tilde{\mathbf{z}}_i^{q,k}\}_{i\in\mathcal{N}}$ be the projected (tangent) points of $\mathbf{z}^{\mathrm{o},k}$ onto each of the lower bound constraints in \eqref{constr:pz_bounds} and \eqref{constr:qz_bounds}, and $\mathcal{Z}^{k}$ the convex restricted set constructed using these tangent points, i.e., $\mathcal{Z}^{k}=\{\mathbf{z}: \eqref{constr:pz&qz}-\eqref{constr:pz_bounds}\}$. To solve \eqref{eqs:convex_opf} with $\mathcal{Z}^{k}$,
we can call a convex solver such as CVXPY~\cite{diamond2016cvxpy} if the objective function is convex, or an NLP solver like IPOPT \cite{wachter2006implementation} if the objective is non-convex. We denote the obtained solution as $\hat{\mathbf{z}}^{k}$, which is surely a feasible point,  and use it to derive a new set of tangent points through projection.
With these updated tangent points, we construct another convex feasible set $\mathcal{Z}^{k+1}$ and solve \eqref{eqs:convex_opf} on it. This iterative process continues until the sequence of objective values converges or a maximum number of steps is reached.
We summarize this iterative algorithm for solving \eqref{eqs:opf} in Table~\ref{tab:algo1}.

\captionsetup[table]{font=small,skip=5pt}
\begin{table}[t!]
\normalsize
\centering
\begin{tabular}{cl}
\hline
\multicolumn{2}{c}{\textbf{Proposed Iterative Algorithm for Solving \eqref{eqs:opf}}}\\
\hline
1:& \textbf{Inputs:}~Initial feasible point $\mathbf{z}^{\mathrm{o},0}$, stopping criterion $\epsilon$, \\
{}& maximum number of iterations $K$.\\
2:& \textbf{For iteration $k$}: \\
3:& Project $\mathbf{z}^{\mathrm{o},k}$ onto $\underline{\mathcal{P}}_i$ and $\underline{\mathcal{Q}}_i$ to get tangents points\\
{ }& for linearization:\\
{ }& \qquad $\tilde{\mathbf{z}}_i^{p,k}=\arg\min_{\mathbf{z}\in\underline{\mathcal{P}}_i} \|\mathbf{z}-\mathbf{z}^{\mathrm{o},k}\|_2^2$ \\
{ }& \qquad $\tilde{\mathbf{z}}_i^{q,k}=\arg\min_{\mathbf{z}\in\underline{\mathcal{Q}}_i} \|\mathbf{z}-\mathbf{z}^{\mathrm{o},k}\|_2^2$. \\
4:& With the obtained tangent points $\{\tilde{\mathbf{z}}_i^{p,k}, \tilde{\mathbf{z}}_i^{q,k}\}_{i\in\mathcal{N}}$, \\
{ }& linearize bus power lower bounds using \eqref{eqn:p_i_lower} and \eqref{eqn:q_i_lower}.\\
5:& Formulate the convex restricted problem \eqref{eqs:convex_opf}.\\
6:& Solve \eqref{eqs:convex_opf} by calling a solver to obtain a solution $\hat{\mathbf{z}}^{k}$\\
{ }& and the associated objective value $\hat{{c}}^{k}$.\\
7:& Update the feasible point by $\hat{\mathbf{z}}^{k}\longrightarrow \mathbf{z}^{\mathrm{o},k+1}$.\\
8:& Let $k+1\longrightarrow k$.\\
9:& \textbf{Repeat the above procedure until $\|\hat{{c}}^{k}-\hat{{c}}^{k-1}\|_2^2\leq \epsilon$}\\
{}& or $k= K$.\\
10:& \textbf{Outputs}: optimal solution $\hat{\mathbf{z}}^{k}$ and objective value  $\hat{{c}}^{k}$.\\
\hline
\end{tabular}
\caption{The iterative algorithm for solving \eqref{eqs:opf} using the proposed convex restriction technique, starting with an initial feasible point $\mathbf{z}^{\mathrm{o}}$.\vspace{-10pt}}
\label{tab:algo1}
\end{table}

\subsection{Convexity of Objective Functions}
In this part, we provide a condition on when the objective functions are convex in the $z$ space, thus making the entire optimization problem in~\eqref{eqs:convex_opf} convex. 

\begin{lemma}\label{lemma:convex_obj}
    If the cost function can be written as 
    \begin{equation} \label{eqn:cost_convex}
    c\big(\mathbf{p}(\mathbf{z}),\mathbf{q}(\mathbf{z})\big) = \sum_i f_i(p_i(\mathbf{z}))+g_i(q_i(\mathbf{z})),
    \end{equation}
    and each $f_i$ and $g_i$ are nondecreasing and convex functions, 
    then the cost $c\big(\mathbf{p}(\mathbf{z}),\mathbf{q}(\mathbf{z})\big)$ is convex in $\mathbf{z}$. 
\end{lemma}
\begin{proof}
    Note that the active and reactive power injections in the $\mathbf{z}$ space, given by \eqref{eqn:active_z} and \eqref{eqn:reactive_z}, are convex functions. Using the fact that composition of a convex function and a convex and nondecreasing function is convex~\cite{boyd2004convex}, we have that the objective function is convex. 
\end{proof}
Lemma \ref{lemma:convex_obj} applies to many standard OPF objective functions that are of the form $c\big(\mathbf{p}(\mathbf{z}),\mathbf{q}(\mathbf{z})\big)=\sum_{i=1}^{N} c_i p_i(\bf{z})$, $c_i\geq 0$, such as minimizing total power loss and total generation cost. It's interesting to note that the conditions on the cost in Lemma~\ref{lemma:convex_obj} are the same as the ones found in SDP~\cite{zhang2013geometry} or SOCP relaxations~\cite{low2014convex}. 

\subsection{Numerical Results}
In this part, we use the iterative algorithm from Table~\ref{tab:algo1} to solve \eqref{eqs:opf} for the IEEE 123-bus distribution network~\cite{Kersting1991RadialDT}, testing with different types of objective functions: loss minimization, generation cost minimization and state estimation. 
We compare our method against several baseline approaches.
The first is running the \textit{runopf} module in MATPOWER~\cite{zimmerman2010matpower}, which is based on Newton-Raphson type algorithms. 
The second is the second-order cone program (SOCP) relaxation~\cite{jabr2006scop,low2014part2}. 
The third method creates a convex region in the power injection space by restricting the permissible voltage phase angles within a polytope and approximating each nonlinear constraint~\cite{lee2019restriction}.
All simulations here are done in Google Colab~\cite{colab} and all codes and data of our experiments are available at \url{https://github.com/zhang-linnng/convex_restriction_transformed}.

\noindent \textbf{Loss Minimization.} The first problem we consider is to minimize the total active power loss in meeting the load. The objective function in this case is 
$$ c\big(\mathbf{p}(\mathbf{z}),\mathbf{q}(\mathbf{z})\big)= \sum_i p_i(\mathbf{z}),$$
which satisfies the condition in Lemma~\ref{lemma:convex_obj}. As a benchmark test, we set the lower bounds in~\eqref{eqs:opf} such that they satisfy the conditions (see~\cite{low2014convex} for details) where the SOCP relaxation is exact. 
Therefore, the SOCP solution is optimal and we are interested in if we can achieve the same loss through convex restriction. 

\begin{figure}[ht]
    \centering
    \includegraphics[scale=0.6]{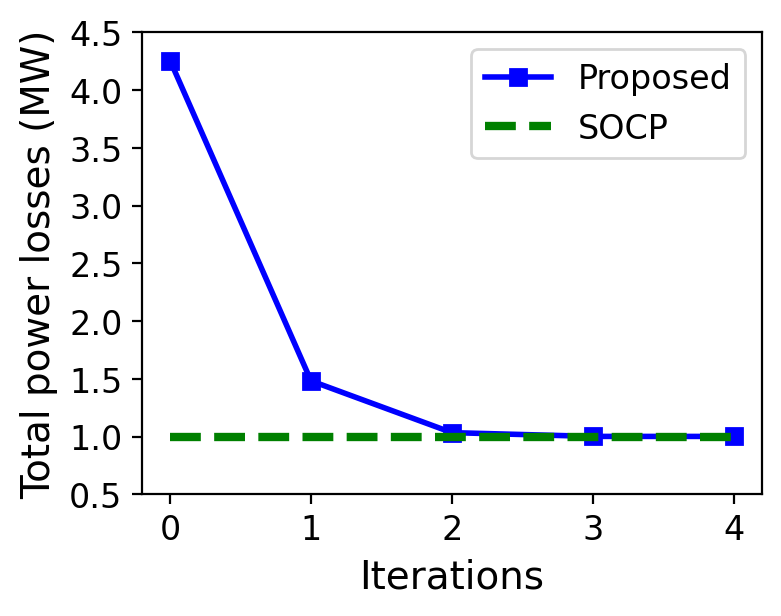}
    \caption{Comparison of the proposed algorithm against the baseline method on the 123-bus network for active power loss minimization.
    The Y-axis is normalized to the minimum objective value within the sequence generated by our algorithm.
    The lower bounds in~\eqref{eqs:opf} are set such that the SOCP relaxation is exact. Our algorithm starts at a higher loss (with a random starting point), but quickly achieves the same loss as SOCP relaxation after 2 iterations.}
    \label{fig:minloss}
\end{figure}
Figure~\ref{fig:minloss} shows the performance of our algorithm compared to the solution of SOCP. We start with a random intialization point that has a high loss, but we quickly reach the same loss as SOCP relaxation in two iterations. This shows that our convex restricted feasible set (after an iteration) contains the optimal solution of the original non-convex problem. It's worth noting that both the SOCP and our algorithm solve convex problems of the same size, and they have roughly the same computational speed. 

\noindent \textbf{Generation Cost Minimization.} Here, we consider a cost in the form of 
$$ c\big(\mathbf{p}(\mathbf{z}),\mathbf{q}(\mathbf{z})\big)= \sum_i c_i p_i(\mathbf{z}),$$
where $c_i$ are positive constants. We set the lower bounds such that the conditions for SOCP relaxation to be exact are not met and test the performance of our algorithm in these scenarios.
Indeed, when the SOCP relaxation is solved, it does not give a physical solution,\footnote{Essentially, the SOCP relaxation relaxes an equality constraint of the type $R_{ij}^2+I_{ij}^2=1$ into $R_{ij}+I_{ij} \leq 1$, where $R_{ij}$ and $I_{ij}$ are variables associated with line $ij$. When the relaxation is exact, the optimal solution $(R_{ij}^{\star})^2$ and $(I_{ij}^{\star})^2$ will satisfy the equality constraint. However, if $(R_{ij}^{\star})^2+(I_{ij}^{\star})^2 <1$, the relaxation is not exact and we cannot recover a solution to the original problem~\eqref{eqs:opf}.} and thus it does not return a feasible solution to the original problem. In addition, the \textit{runopf} module in MATPOWER also fails to find a solution. As for the method that assumes a polytopic phase angle feasible set, it proves to be overly restrictive and results in an empty set unless we carefully adjust the hyperparameters used for approximating each nonlinear constraint.

\begin{figure}[ht]
    \centering
    \includegraphics[scale=0.6]{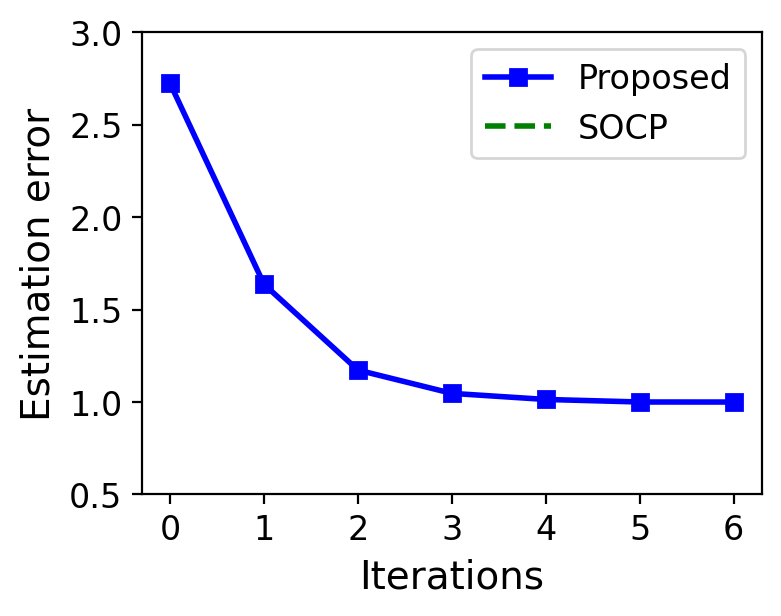}
    \caption{Performance of the proposed algorithm on the 123-bus network for generation cost minimization. 
    The Y-axis is normalized to the minimum objective value within the sequence generated by our algorithm.
    SOCP turns out to be inexact in this case and does not produce a feasible solution. 
    The \textit{runopf} routine in MATPOWER also fails to converge. 
    In contrast, our algorithm reliably decreases the cost and converges to a physically feasible solution.
    }
    \label{fig:mincost}
\end{figure}
In Figure~\ref{fig:mincost}, we show the performance of our algorithm, where
our method reliably decreases the cost and converges to a good feasible solution. 
This shows the benefit of working with a convex restriction, since we are always guaranteed to find a feasible solution of the original problem.

\noindent \textbf{State Estimation}
Given some measurements of nodal active and reactive power injections, the goal of state estimation is to find the optimal $\bf{z}$ that best matches these measurements. For example, suppose we take active and reactive power measurements at each bus, denoted by $\{\hat{p}_1, \ldots, \hat{p}_{N}, \hat{q}_1, \ldots, \hat{q}_{N}\}$, state estimation solves \eqref{eqs:convex_opf}  with the following objective:
\begin{equation} \label{eqn:state_estimation}
c\big(\mathbf{p}(\mathbf{z}),\mathbf{q}(\mathbf{z})\big) = \sum_{i=1}^{N} (\hat{p}_i-p(\mathbf{z})_i)^2 + \sum_{i=1}^{N} (\hat{q}_i-q(\mathbf{z})_i)^2.
\end{equation}
Due to the non-monotonicity in the objective function \eqref{eqn:state_estimation}, \eqref{eqs:convex_opf} is non-convex, but it can still be solved using a nonlinear solver such as IPOPT. Note that since the feasible set of the convex restricted problems are convex, we are still guaranteed that the solution is feasible
(although perhaps not optimal). 

\begin{figure}[ht]
    \centering
    \includegraphics[scale=0.6]{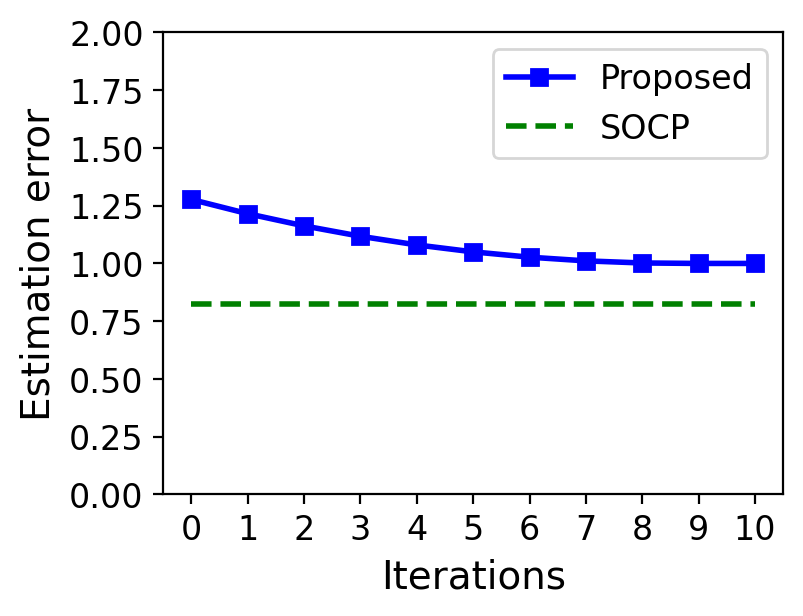}
    \caption{Performance of the proposed algorithm on the 123-bus network for state estimation. 
    The Y-axis is normalized to the minimum objective value within the sequence generated by our algorithm.
    SOCP turns out to be inexact and the cost returned by it is not actually achievable.
    But it serves as a lower bound, and our algorithm is not very far from this bound, and it produces a sequence of feasible solutions. 
    }
    \label{fig:state_estimation}
\end{figure}
In the case of state estimation, the SOCP relaxation is quite far away from being exact. 
This is expected, since the objective function is not increasing in the active power, and convex relaxation algorithms tend to struggle to produce physically meaningful solutions. In contrast, our proposed algorithm continues to perform effectively. 
To illustrate this, we present the convergence of the objective value sequence generated by our algorithm in Figure \ref{fig:state_estimation}. 
Since the SOCP is inexact, the dashed line in Figure~\ref{fig:state_estimation} is not actually achievable. But it serves as a lower bound, and our algorithm converges to a physically feasible solution not very far from this bound.

\section{Conclusions}\label{sec:conclusions}
In this paper, we focused on developing a convex restriction approach for solving AC power flow problems in radial networks.
We introduced a simple change of variables technique, showing that the active and reactive power equations are naturally convex in the transformed coordinate space. 
A detailed procedure was provided to construct a convex subset in this transformed space, and the convex restriction constructed in this way was shown to be a maximal one.
Furthermore, we proposed an iterative algorithm to improve on the solution quality by constructing a series of convex restricted sets, each containing solutions progressively closer to the optimal ones. 
We conducted numerical experiments on the IEEE 123-bus distribution network and solved the ACOPF problem by applying different types of objective functions. 
The numerical results showed that our method produced good feasible solutions within just a few iterations for all study cases, even when traditional methods failed to return a solution.

\bibliographystyle{IEEEtran}
\bibliography{refs.bib}






\end{document}